\documentclass[a4paper, 12pt]{article}
\usepackage{amssymb,latexsym}
\usepackage[cp1251]{inputenc}
\usepackage[english]{babel}
\usepackage{fullpage,amsfonts,amsthm}
\usepackage[intlimits]{amsmath}
\usepackage{indentfirst}
\newtheorem{theorem}{Theorem}
\newtheorem{lemma}[theorem]{Lemma}

\title{
\vspace{15mm}
\textbf{On an equation involving fractional powers with prime numbers of a special type}}
\author{
\vspace{8mm}
Zh. H. Petrov \footnote {Supported by SU Grant 80-10-215/2017. }}
\date{}

\begin{document}

\maketitle
\begin{abstract}
We consider the equation $[p_{1}^{c}] + [p_{2}^{c}] + [p_{3}^{c}] = N$, where $N$ is a sufficiently large integer,
and prove that if $1 < c < \frac{17}{16}$, then it has a solution in prime numbers $p_{1}$, $p_{2}$, $p_{3}$
such that each of the numbers $p_{1} + 2$, $p_{2} + 2$, $p_{3} + 2$ has at most $\left [ \frac{95}{17 - 16c} \right ]$
prime factors, counted with the multiplicity.
\end{abstract}

\section{Introduction and statement of the result}

In 1937 I. M. Vinogradov \cite{4} proved that for every sufficiently large odd integer $N$ the equation
\begin{equation} \label{0.01}
p_{1} + p_{2} + p_{3} = N
\end{equation}
has a solution in prime numbers $p_{1}$, $p_{2}$, $p_{3}$.

Analogous problem involving diophantine inequality was considered in 1952 by Piatetski-Shapiro \cite{20}. 
In 1992, Tolev~\cite{12} established that
if $1 < c < \frac{15}{14}$, then the diophantine inequality 
\begin{equation*}
|p_{1}^{2} + p_{2}^{c} + p_{3}^{c} - N| < N^{-\kappa }
\end{equation*}
has a solution in prime numbers $p_{1}$, $p_{2}$, $p_{3}$ for certain $\kappa = \kappa(c) > 0$. Several improvements
were made and the strongest of them is due to Baker and Weingartner \cite{13}, who improved Tolev's result with
$1 < c < \frac{10}{9}$.

In 1995, M. B. Laporta and D. I. Tolev \cite{5} considered the equation
\begin{equation} \label{0.1}
[p_{1}^{c}] + [p_{2}^{c}] + [p_{3}^{c}] = N,
\end{equation}
where $c \in \mathbb{R}$, $c > 1$, $N \in \mathbb{N}$ and $[t]$ denotes the integer part of $t$.
They showed that if $1 < c < \frac{17}{16}$ and
$N$ is a sufficiently large integer, then the equation \eqref{0.1} has a solution in prime numbers
$p_{1}$, $p_{2}$, $p_{3}$.

For any natural number $r$, let $\mathcal P_{r}$ denote the set of $r$-almost primes, i.e.
the set of natural numbers having at most $r$ prime factors counted with multiplicity.
There are many papers devoted to the study of
problems involving primes and almost primes. For example, in 1973 J. R. Chen \cite{6} established that there exist infinitely many primes $p$
such that $p + 2 \in \mathcal P_{2}$.
In 2000 Tolev \cite{7} proved that for every sufficiently large integer $N~\equiv~3~\pmod 6$ the equation \eqref{0.01}
has a solution in prime numbers $p_{1}$, $p_{2}$, $p_{3}$ such that
$p_{1} + 2 \in P_{2}$, $p_{2} + 2 \in P_{5}$, $p_{3} + 2 \in P_{7}$. Thereafter this result was improved by
Matom\"{a}ki and Shao~\cite{8}, who showed that for every sufficiently large integer $N \equiv 3 \pmod 6$ the equation \eqref{0.01}
has a solution in prime numbers $p_{1}$, $p_{2}$, $p_{3}$ such that
$p_{1} + 2, p_{2} + 2, p_{3} + 2 \in \mathcal P_{2}$.


Recently Tolev \cite{3} established that if $N$ is sufficiently large, $E > 0$ is an arbitrarily
large constant and $1 < c < \frac{15}{14}$ then the inequality
\begin{equation*}
|p_{1}^{2} + p_{2}^{c} + p_{3}^{c} - N| < (\log N)^{-E}
\end{equation*}
has a solution in prime numbers $p_{1}$, $p_{2}$, $p_{3}$, such that each of the numbers
$p_{1} + 2$, $p_{2} + 2$, $p_{3} + 2$ has at most $\left [ \frac{369}{180 - 168c} \right ]$
prime factors, counted with the multiplicity.

In this paper, we prove the following
\begin{theorem}\label{theorem 1}
Suppose that $1 < c < \frac{17}{16}$. Then for every sufficiently large $N$ the equation \eqref{0.1}
has a solution in prime numbers $p_{1}$, $p_{2}$, $p_{3}$,
such that each of the numbers $p_{1} + 2$, $p_{2} + 2$, $p_{3} + 2$ has at most $\left [ \frac{95}{17 - 16c} \right ]$
prime factors, counted with the multiplicity.
\end{theorem}

\vspace{5mm}
We note that the integer $\left [ \frac{95}{17 - 16c} \right ]$ is equal to $95$ if $c$ is close to $1$ and it is
large if $c$ is close to $\frac{17}{16}$.

To prove Theorem \ref{theorem 1} we combine ideas developed by Laporta and Tolev \cite{5} and Tolev~\cite{3}.
First we apply a version of the vector sieve and then the
circle method. In section~4 we find an asymptotic formula for the integrals $\Gamma _{1}'$ and $\Gamma _{4}'$ (defined by
\eqref{3.7} and \eqref{3.9} respectively). In section~5 we estimate $\Gamma _{1}''$ and $\Gamma _{4}''$ (defined by
\eqref{3.8} and \eqref{3.10} respectively) and we will complete the proof.

\vspace{5mm}
\textbf{Acknowledgements.} The author wishes to express his thanks to Professor D. Tolev 
for suggesting the problem and for the helpful conversations. 

\section{Notations and some lemmas}
We use the following notations: with $[t]$ we denote the integer part of $t$ and $\{t\} = t - [t]$ is the
fractional part of $t$. With $||t||$ we denote the distance from $t$ to the nearest integer. As usual with
$\mu (n)$, $\varphi (n)$ and $\Lambda (n)$ we denote respectively, M\"{o}bius' function, Euler's function
and von Mangoldt's function. Also $e(t) = e^{2\pi i t}$.

We use Vinogradov's notation $A \ll B$, which is equivalent to $A = O(B)$. If we have
simultaneously $A \ll B$ and $B \ll A$, then we shall write $A \asymp B$.

We reserve $p, p_{1}, p_{2}, p_{3} $ for prime numbers. By $\epsilon $ we denote an arbitrarily small
positive number, which is not necessarily the same in the different formulae.

With $\mathbb{N}$, $\mathbb{Z}$ and $\mathbb{R}$ we will denote respectively the set of natural numbers, the set of
integer numbers and the set of real numbers.

Now we introduce some lemmas, which shall be used later.

\vspace{5mm}
\begin{lemma} \label{lemma 1}
Suppose that $D \in \mathbb{R}, D > 4$. There exist arithmetical functions $\lambda ^{\pm}(d)$ (called Rosser's functions of level D)
with the following properties:
\begin{enumerate}
\item
    For any positive integer $d$ we have
    \begin{equation*}
    |\lambda ^{\pm}(d)| \leq 1, \qquad \qquad \lambda ^{\pm}(d) = 0 \quad \text{if} \quad d > D \quad \text{or} \quad \mu (d) = 0.
    \end{equation*}
\item
    If $n \in \mathbb{N}$ then
    \begin{equation*}
    \sum_{d|n}{\lambda^{-}(d)} \leq \sum_{d|n}{\mu (d)} \leq \sum_{d|n}{\lambda^{+}(d)}.
    \end{equation*}
\item
    If $z \in \mathbb{R}$ is such that $z^{2} \leq D \leq z^{3}$ and if
    \begin{equation} \label{1.01}
    P(z) = \prod_{2 < p < z}{p}, \; \; \mathcal {B} = \prod_{2 < p < z}{\left ( 1 - \frac{1}{p - 1} \right )}, \; \;
	\mathcal {N}^{\pm } = \sum_{d|P(z)}{\frac{\lambda ^{\pm} (d)}{\varphi (d)}}, \; \; s_{0} = \frac{\log D}{\log z},
    \end{equation}
then we have
\begin{align}
\mathcal {B} & \leq \mathcal {N} ^{+} \leq
				\mathcal {B} \left ( F(s_{0}) + O\left ( (\log D)^{-\frac{1}{3}} \right )\right ), \label{1.02}\\
\mathcal {B} & \geq \mathcal {N} ^{-} \geq
				\mathcal {B} \left ( f(s_{0}) + O\left ( (\log D)^{-\frac{1}{3}} \right )\right ) \label{1.03},
\end{align}
where $F(s)$ and $f(s)$ satisfy
\begin{equation} \label{1.04}
f(s) = 2e^{\gamma }s^{-1}\log (s - 1), \quad F(s) = 2e^{\gamma }s^{-1} \quad \text{ for } \quad 2 \leq s \leq 3.
\end{equation}
Here $\gamma $ is Euler's constant.
\end{enumerate}
\end{lemma}
\textbf{Proof.} See Greaves \cite[Chapter 4]{2}. \hfill $ \square $ \\

\begin{lemma} \label{lemma 2}
Suppose that $\Lambda_{i}, \Lambda_{i}^{\pm}$ are real numbers satisfying $\Lambda_{i} = 0$ or $1$,
$\Lambda_{i}^{-} \leq \Lambda_{i} \leq \Lambda_{i}^{+}$, $i = 1, 2, 3$. Then
\begin{equation} \label{1.1}
\Lambda_{1} \Lambda _{2} \Lambda _{3} \geq \Lambda_{1}^{-}\Lambda_{2}^{+}\Lambda_{3}^{+} +
                                            \Lambda_{1}^{+}\Lambda_{2}^{-}\Lambda_{3}^{+} +
                                            \Lambda_{1}^{-}\Lambda_{2}^{+}\Lambda_{3}^{-} -
                                            2\Lambda_{1}^{+}\Lambda_{2}^{+}\Lambda_{3}^{+}.
\end{equation}
\end{lemma}
\textbf{Proof.} The proof is similar to the proof of Lemma 2 in \cite{1}. \hfill $ \square $ \\

\begin{lemma} \label{lemma 3}
Suppose that $x$, $y \in \mathbb{R}$ and $M \in \mathbb{N}$, $M \geq 3$. Then
\begin{equation*}
e(-x\{ y \}) = \sum_{|m| \leq M}{c_{m}e(my)} + O\left ( \min \left ( 1, \frac{1}{M||y||}\right )\right ),
\end{equation*}
where
\begin{equation} \label{100}
c_{m} = \frac{1 - e(-x)}{2 \pi i (x + m)}.
\end{equation}
\end{lemma}
\textbf{Proof.} Proof can be find in Buriev \cite[Lemma 12]{16}. \hfill $ \square $ \\

\begin{lemma} \label{lemma 5}
Consider the integral
\begin{equation*}
I = \int_{a}^{b}e(f(x))dx,
\end{equation*}
where $f(x)$ is real function with continuous second derivative and monotonous first derivative. If
$|f'(x)| \geq~h~>~0$, for all $x \in [a, b]$ then $I \ll h^{-1}$.
\end{lemma}
\textbf{Proof.} See  \cite[p. 71]{15}. \hfill $ \square $ \\

\section{Beginning of the proof}

Let $\eta $, $\delta $, $\xi $ and $\mu $ are positive real numbers depending on $c$. We shall specify them later. Now we only
assume that they satisfy the conditions
\begin{equation} \label{3.1}
\xi + 3\delta < \frac{12}{25}, \qquad 2 < \frac{\delta }{\eta } < 3, \qquad \mu < 1.
\end{equation}
We denote
\begin{equation} \label{3.2}
X = N^{\frac{1}{c}}, \qquad z = X^{\eta }, \qquad D = X^{\delta }, \qquad \Delta = X^{\xi - c}
\end{equation}
and
\begin{equation} \label{3.3}
P(z) = \prod_{2 < p < z}{p}.
\end{equation}

Consider the sum
\begin{equation} \label{3.4}
\Gamma = \sum_{\substack{\mu X < p_{1}, p_{2}, p_{3} \leq X\\
						\\
						 [p_{1}^{c}] + [p_{2}^{c}] + [p_{3}^{c}] = N \\
						 \\
						 (p_{i} + 2, P(z)) = 1, \; i = 1, 2, 3}}
{(\log p_{1})(\log p_{2})(\log p_{3})}.
\end{equation}

If we prove the inequality
\begin{equation} \label{3.5}
\Gamma > 0,
\end{equation}
then the equation \eqref{0.1} would have a solution in primes $p_{1}, p_{2}, p_{3}$ satisfying conditions in the sum $\Gamma $.
Suppose that $p_{i} + 2$ has $l$ prime factors, counted with multiplicity. From \eqref{3.2}, \eqref{3.3} and
$(p_{i} + 2, P(z)) = 1$ we have
\begin{equation*}
X + 2 \geq p_{i} + 2 \geq z^{l} = X^{\eta l}
\end{equation*}
and then $l < \frac{1}{\eta }$. This means that $p_{i} + 2$ has at most
$[\eta ^{-1}]$ prime factors counted with multiplicity. Therefore, to prove Theorem \ref{theorem 1} we have to establish \eqref{3.5} for an appropriate
choice of $\eta $.

For $i = 1, 2, 3$ we define
\begin{equation} \label{2.1}
\Lambda _{i} = \sum_{d|(p_{i} + 2, P(z))}{\mu(d)} =
		\begin{cases}
		1 & \quad \text{ if } (p_{i} + 2, P(z)) = 1, \\
		0 & \quad \text{ otherwise. }
		\end{cases}
\end{equation}
Then we find that
\begin{equation*}
\Gamma = \sum_{\substack{\mu X < p_{1}, p_{2}, p_{3} \leq X \\ [p_{1}^{c}] + [p_{2}^{c}] + [p_{3}^{c}] = N}}
{\Lambda _{1}\Lambda _{2} \Lambda _{3}(\log p_{1})(\log {p_{2}})(\log {p_{3}})}.
\end{equation*}
We can write $\Gamma $ as
\begin{equation*}
\Gamma = \sum_{\mu X < p_{1}, p_{2}, p_{3} \leq X}
{\Lambda _{1}\Lambda _{2} \Lambda _{3}(\log p_{1})(\log {p_{2}})(\log {p_{3}})
\int_{-\frac{1}{2}}^{\frac{1}{2}}}{e(\alpha ([p_{1}^{c}] + [p_{2}^{c}] + [p_{3}^{c}] - N))d\alpha }.
\end{equation*}

Suppose that $\lambda ^{\pm }(d)$ are the Rosser functions of level $D$. Let also denote
\begin{equation} \label{2.2}
\Lambda _{i} ^{\pm } = \sum_{d|(p_{i} + 2, P(z))}{\lambda ^{\pm } (d)}, \hspace{10mm} i = 1, 2, 3.
\end{equation}
Then from Lemma \ref{lemma 1}, \eqref{2.1} and \eqref{2.2} we find that
\begin{equation*}
\Lambda _{i}^{-} \leq \Lambda _{i} \leq \Lambda _{i}^{+}.
\end{equation*}

We use Lemma \ref{lemma 2} and find that
\begin{equation*}
\Gamma \geq \Gamma _{1} + \Gamma _{2} + \Gamma _{3} - 2\Gamma _{4},
\end{equation*}
where $\Gamma _{1}, \dots , \Gamma _{4}$ are the contributions coming from the consecutive terms of the right side of~\eqref{1.1}.
We have $\Gamma_{1} = \Gamma_{2} = \Gamma_{3} $ and
\begin{align*}
\Gamma_{1}
	& = \sum_{\mu X < p_{1}, p_{2}, p_{3} \leq X}
		{\Lambda _{1}^{-}\Lambda _{2}^{+}\Lambda _{3}^{+}(\log p_{1})(\log p_{2})(\log p_{3})
		\int_{-\frac{1}{2}}^{\frac{1}{2}}{e(\alpha([p_{1}^{c}] + [p_{2}^{c}] + [p_{3}^{c}] - N))d\alpha}}, \\
\Gamma_{4}
	& = \sum_{\mu X < p_{1}, p_{2}, p_{3} \leq X}
		{\Lambda _{1}^{+}\Lambda _{2}^{+}\Lambda _{3}^{+}(\log p_{1})(\log p_{2})(\log p_{3})
		\int_{-\frac{1}{2}}^{\frac{1}{2}}{e(\alpha([p_{1}^{c}] + [p_{2}^{c}] + [p_{3}^{c}] - N))d\alpha}}.
\end{align*}
Hence, we get
\begin{equation} \label{3.5.1}
\Gamma \geq 3\Gamma_{1} - 2\Gamma_{4}.
\end{equation}

Let first consider $\Gamma_{1}$. We have
\begin{equation} \label{3.6}
\Gamma_{1} = \int_{-\frac{1}{2}}^{\frac{1}{2}}{e(-N\alpha )L^{-}(\alpha )L^{+}(\alpha )^{2}d\alpha },
\end{equation}
where
\begin{equation*}
L^{\pm}(\alpha ) = \sum_{\mu X < p \leq X}{(\log p)e(\alpha [p^{c}])}\sum_{d|(p+2, P(z))}{\lambda ^{\pm}(d)}.
\end{equation*}
Changing the order of summation in $L^{\pm}(\alpha )$, we get
\begin{equation*}
L^{\pm}(\alpha ) = \sum_{d|P(z)}{\lambda ^{\pm}(d)}
\sum_{\substack{\mu X < p \leq X \\ p + 2 \equiv 0 (\bmod d)}}{(\log p )e(\alpha [p^{c}])}.
\end{equation*}

We divide the integral from \eqref{3.6} into two parts:
\begin{equation} \label{3.5.2}
\Gamma _{1} = \Gamma _{1}' + \Gamma _{1} '',
\end{equation}
where
\begin{align}
\Gamma _{1}'
			& = \int_{|\alpha | < \Delta }{e(-N\alpha )L^{-}(\alpha )L^{+}(\alpha )^{2}d\alpha }, \label{3.7}\\
\Gamma _{1}''
			& = \int_{\Delta < |\alpha | < \frac{1}{2}}{e(-N\alpha )L^{-}(\alpha )L^{+}(\alpha )^{2}d\alpha }. \label{3.8}
\end{align}

Similarly, for $\Gamma _{4}$ we have
\begin{equation} \label{3.5.3}
\Gamma _{4} = \Gamma _{4}' + \Gamma _{4} '',
\end{equation}
where
\begin{align}
\Gamma _{4}'
			& = \int_{|\alpha | < \Delta }{e(-N\alpha )L^{+}(\alpha )^{3}d\alpha }, \label{3.9}\\
\Gamma _{4}''
			& = \int_{\Delta < |\alpha | < \frac{1}{2}}{e(-N\alpha )L^{+}(\alpha )^{3}d\alpha }. \label{3.10}
\end{align}

\section{The integrals $\Gamma _{1}'$ and $\Gamma _{4}'$}

We shall find an asymptotic formula for the integrals $\Gamma _{1}'$ and $\Gamma _{4}'$ defined by
\eqref{3.7} and \eqref{3.9}, respectively. The arithmetic structure of the Rosser weights $\lambda ^{\pm }(d)$ are not important here, so we consider a sum of the form
\begin{equation} \label{4.1}
L(\alpha ) = \sum_{d \leq D}{\lambda (d)}
			 \sum_{\substack{\mu X < p \leq X \\ p + 2 \equiv 0 (\bmod d)}}{(\log p)e(\alpha [p^{c}])},
\end{equation}
where $\lambda (d)$ are real numbers satisfying
\begin{equation} \label{4.100}
|\lambda(d)| \leq 1, \qquad \lambda (d) = 0 \quad \text{ if } \quad 2|d \quad \text{ or } \quad \mu (d) = 0.
\end{equation}
It is easy to see that
\begin{align*}
L(\alpha )
			& = \sum_{d \leq D}{\lambda (d)}
				\sum_{\substack{\mu X < p \leq X \\ p + 2 \equiv 0 (\bmod d)}}{(\log p)e(\alpha p^{c} + O(|\alpha |))} \\
			& = \sum_{d \leq D}{\lambda (d)}
				\sum_{\substack{\mu X < p \leq X \\ p + 2 \equiv 0 (\bmod d)}}{(\log p)e(\alpha p^{c})(1 + O(|\alpha |))} \\
			& = \sum_{d \leq D}{\lambda (d)}
				\sum_{\substack{\mu X < p \leq X \\ p + 2 \equiv 0 (\bmod d)}}{(\log p)e(\alpha p^{c})}
				+ O\left ( \sum_{d \leq D}{\sum_{\substack{p \leq X \\ p + 2 \equiv 0 (\bmod d)}}{(\log p )|\alpha |}}\right ).
\end{align*}
If $D \leq X$ then
\begin{equation} \label{4.2}
L(\alpha ) = \overline{L}(\alpha ) + O(\Delta X (\log X)),
\end{equation}
where
\begin{equation*}
\overline{L} = \sum_{d \leq D}{\lambda (d)}
				\sum_{\substack{\mu X < p \leq X \\ p + 2 \equiv 0 (\bmod d)}}{(\log p)e(\alpha p^{c})}.
\end{equation*}

For $\overline{L}(\alpha )$ we use the asymptotic formula from Lemma 10 in \cite{3}. From \eqref{3.1} and
\eqref{3.2} we see that, when $|\alpha| < \Delta$, then for every constant $A > 0$, we have
\begin{equation} \label{4.3}
\overline{L}(\alpha ) = \sum_{d \leq D}{\frac{\lambda (d)}{\varphi (d)}I(\alpha )} + O(X(\log X)^{-A}),
\end{equation}
where
\begin{equation} \label{4.3.5}
I(x) = \int_{\mu X}^{X}{e(\alpha t^{c})dt}.
\end{equation}
Hence from \eqref{3.2}, \eqref{4.2} and \eqref{4.3} we see that if $\xi < c$ then
\begin{equation} \label{4.4}
L(\alpha )  = \sum_{d \leq D}{\frac{\lambda (d)}{\varphi (d)}I(\alpha )} + O(X(\log X)^{-A}).
\end{equation}

If $|\alpha | < \Delta $, then from \eqref{1.01} and \eqref{4.4} we find
\begin{equation} \label{4.5}
L^{\pm }(\alpha ) = \mathcal {N}^{\pm } I(\alpha ) + O(X(\log X)^{-A}).
\end{equation}
Let
\begin{equation} \label{4.6}
\mathcal {M}^{\pm } = \mathcal {M}^{\pm }(\alpha ) = \sum_{d \leq D}{\frac{\lambda ^{\pm }(d)}{\varphi (d)}I(\alpha )}
= \mathcal {N}^{\pm }I(\alpha ).
\end{equation}
It is easy to see that
\begin{equation} \label{4.6.5}
\mathcal {M}^{\pm } \ll (\log X)|I(\alpha )|.
\end{equation}

We use \eqref{4.5}, \eqref{4.6} and the identity
\begin{equation*}
L^{-}(L^{+})^{2} = (L^{-} - \mathcal {M}^{-})(L^{+})^{2} + (L^{+} - \mathcal {M}^{+})\mathcal {M}^{-}L^{+}
					+ (L^{+} - \mathcal {M}^{+})\mathcal {M}^{+}\mathcal {M}^{-} + \mathcal {M}^{-}(\mathcal {M}^{+})^{2},
\end{equation*}
to find that
\begin{equation} \label{4.7}
|L^{-}(L^{+})^{2} - \mathcal {M}^{-}(\mathcal {M}^{+})^{2}| \ll
X(\log X)^{-A}\left ( |L^{+}|^{2} + |\mathcal {M}^{-}|^{2} + |\mathcal {M}^{+}|^{2} \right ).
\end{equation}

Let
\begin{equation} \label{4.8}
B = \int_{|\alpha | < \Delta }{e(-N\alpha )\mathcal M ^{-}(\alpha ) (\mathcal M^{+}(\alpha ))^{2}d\alpha }.
\end{equation}
From \eqref{3.7}, \eqref{4.6.5} -- \eqref{4.8} we have
\begin{equation*}
\Gamma _{1}' - B \ll
X(\log X)^{2 - A}\left ( \int_{|\alpha | < \Delta }{|L^{+}(\alpha )|^{2}d\alpha }
+ \int_{|\alpha| < \Delta }{|I(\alpha )|^{2}d\alpha }\right ).
\end{equation*}

We need the next lemma, which is an analog of Lemma 11 in \cite{3}.
\begin{lemma} \label{lemma 4.1}
If $\Delta \leq X^{1 - c}$, then for the sum $L(\alpha )$ defined by \eqref{4.1} and for the integral $I(\alpha )$
defined by \eqref{4.3.5} we have
\begin{align*}
\int_{|\alpha | < \Delta }{|L(\alpha )|^{2}d\alpha }
		& \ll X^{2 - c}(\log X)^{6}, \\
\int_{|\alpha | < \Delta }{|I(\alpha )|^{2}d\alpha }
		& \ll X^{2 - c}(\log X)^{6}, \\
\int_{|\alpha | < 1 }{|I(\alpha )|^{2}d\alpha }
		& \ll X(\log X)^{5}.
\end{align*}
\end{lemma}
\textbf{Proof.} The proof is similar to the proof of Lemma 11 in \cite{3}.   \hfill $ \square $ \\

Hence
\begin{equation} \label{3.5.99}
\Gamma _{1}' - B \ll X^{3 - c}(\log X)^{8 - A}.
\end{equation}

Consider now the integral
\begin{equation} \label{3.5.100}
B_{1} = \int_{-\infty }^{\infty }{e(-N\alpha )I(\alpha )^{3}d\alpha }.
\end{equation}
Using the method in Lemma 5.6.1 in \cite{14} we find
\begin{equation} \label{3.5.7}
B_{1} \gg X^{3 - c}.
\end{equation}

For $I(\alpha )$ we apply Lemma \ref{lemma 5} and see that $I(\alpha ) \ll |\alpha |^{-1}X^{1 - c}$.
Then from \eqref{3.2}, \eqref{4.6}, \eqref{4.8} and \eqref{3.5.100} we find
\begin{equation} \label{3.5.101}
|\mathcal {N}^{-}(\mathcal {N}^{+})^{2}B_{1} - B|
\ll (\log X)^{3}\int_{|\alpha | > \Delta }{|I(\alpha )|^{3}d\alpha } \ll (\log x)^{3} X^{3 - c - 2\xi}.
\end{equation}
If $A = 12$, then using \eqref{3.5.99} and \eqref{3.5.101} we find
\begin{equation} \label{3.5.5}
\Gamma _{1}' = \mathcal {N}^{-}(\mathcal {N}^{+})^{2}B_{1} + O(X^{3 - c}(\log X)^{- 4}).
\end{equation}
We proceed with $\Gamma _{2}'$ in the same way and prove that
\begin{equation} \label{3.5.6}
\Gamma _{4}' = (\mathcal {N}^{+})^{3}B_{1} + O(X^{3 - c}(\log X)^{- 4}).
\end{equation}

\section{The estimation of the integrals $\Gamma _{1}''$ and $\Gamma _{4}''$ and the end of the proof}
In this section we consider the integrals $\Gamma _{1}''$ and $\Gamma _{4}''$ defined by
\eqref{3.8} and \eqref{3.10} respectively. We show that the integrals $\Gamma _{1}''$ and $\Gamma _{4}''$
are small enough. Now we assume that
\begin{equation} \label{5.1}
	\xi = \frac{16c - 5}{32}, \qquad \qquad \delta = \frac{17 - 16c}{32}.
\end{equation}
It is obvious that for $\Gamma _{1}''$ defined by \eqref{3.8} we have
\begin{equation*}
\Gamma _{1}'' \ll \max_{\Delta \leq |\alpha | \leq \frac{1}{2}}{|L^{-}(\alpha )|}
\int_{0}^{1}{|L^{+}(\alpha )|^{2}d\alpha }.
\end{equation*}
We use Lemma \ref{lemma 4.1} and find that
\begin{equation} \label{5.1.0}
\Gamma _{1}'' \ll X(\log X)^{5}\max_{\Delta \leq |\alpha | \leq \frac{1}{2}}{|L^{-}(\alpha )|}.
\end{equation}
From \eqref{4.1} we see that
\begin{equation} \label{5.1.1}
L(\alpha ) = L_{1}(\alpha ) + O\left ( X^{\frac{1}{2} + \varepsilon} \right ),
\end{equation}
where
\begin{equation*}
L_{1}(\alpha ) = \sum_{d \leq D}{\lambda (d)}\sum_{\substack{\mu X < n \leq X \\ n + 2 \equiv 0(\bmod d)}}
{\Lambda (n)e(\alpha [n^{c}])}.
\end{equation*}

Let $M = X^{\kappa }$, for some $\kappa$, which will be specified later.
Now for $L_{1}(\alpha )$ applying Lemma~\ref{lemma 3} with parameters $x = \alpha $, $y = n^{c}$ and $M$
(we note that $[t] = t - \{ t \}$). Hence
\begin{align}
L_{1}(\alpha )
				= \sum_{|m| \leq M}{c_{m}}\sum_{d \leq D}{\lambda (d)}
				  \sum_{\substack {\mu X < n \leq X \\ n + 2 \equiv 0 (\bmod d)}}{\Lambda (n)e((\alpha + m)n^{c})} + \notag \\
				  + O\left ( X^{\varepsilon }\sum_{\mu X < n \leq X}{\min \left ( 1, \frac{1}{M||n^{c}||}\right )} \right ). \label{5.2}
\end{align}
We need the following
\begin{lemma} \label{lemma 10}
Suppose that $D$, $\Delta $ are define by \eqref{3.2} and $\xi $, $\delta $ are specified by \eqref{5.1}.
Suppose also that $\lambda (d)$ satisfy \eqref{4.100} and $c_{m}$ are define by \eqref{100}. Then
\begin{align*}
& \max _{\Delta \leq \alpha \leq M + 1}
			 \left |  \sum_{|m| \leq M}{c_{m}}\sum_{d \leq D}{\lambda (d)}
			  \sum_{\substack {\mu X < n \leq X \\ n + 2 \equiv 0(\bmod d)}}{\Lambda (n)e(\alpha n^{c})} \right | \ll \\
			& \ll x^{\varepsilon }
			\left ( X^{\frac{1}{3} + \frac{c}{2}}DM^{\frac{1}{2}}
			+ X^{1 - \frac{c}{2}}\Delta ^{-\frac{1}{2}}
			+ X^{\frac{3}{4} + \frac{c}{6}}D^{\frac{2}{3}}M^{\frac{1}{6}}
			+ X^{\frac{5}{6}}
			+ X^{1 -\frac{c}{6}}D^{\frac{1}{3}}\Delta ^{-\frac{1}{6}}
			+ X^{1 - \frac{c}{4}}\Delta ^{-\frac{1}{4}}
			\right ).
\end{align*}				
\end{lemma}
\textbf{Proof.} See Lemma 15 in \cite{3}. \hfill $ \square $ \\
Also we need the estimate
\begin{equation} \label{5.3}
\sum_{\mu X < n \leq X}{\min \left ( 1, \frac{1}{M||n^{c}||}\right )}
			\ll X^{\varepsilon }\left ( XM^{-1} + M^{\frac{1}{2}}X^{\frac{c}{2}} \right ).
\end{equation}
For the proof we note that the Fourier series of $\min \left ( 1, \frac{1}{M||n^{c}||}\right )$ is given by
\begin{equation} \label{1000}
\min \left ( 1, \frac{1}{M||n^{c}||}\right ) = \sum_{k \in \mathbb{N}}{b_{M}(k)e(kn^{c})},
\end{equation}
where the Fourier coefficients satisfy
\begin{equation} \label{1002}
|b_{M}(k)| \leq
			\begin{cases}
			\frac{4\log M}{M} & \quad \text{ if } k \in \mathbb{Z}, \\
			\frac{M}{n^{2}} & \quad \text{ if } k \in \mathbb{Z}, n \neq 0.
			\end{cases}
\end{equation}
From \eqref{1000} we get
\begin{equation}\label{1001}
\sum_{\mu X < n \leq X}{\min \left ( 1, \frac{1}{M||n^{c}||}\right )} =
		\sum_{\mu X < n \leq X}{\sum_{k \in \mathbb{N}}{b_{M}(k)e(kn^{c})}}.
\end{equation}
Changing the order of summation in last formula we get
\begin{equation*}
\sum_{\mu X < n \leq X}{\min \left ( 1, \frac{1}{M||n^{c}||}\right )} =
		\sum_{k \in \mathbb{N}}{b_{M}(k)H(k)},
\end{equation*}
where
\begin{equation*}
H(k) = \sum_{\mu X < n \leq X}{e(kn^{c})}.
\end{equation*}
Now using \eqref{1002} and \eqref{1001} and the identity $|H(k)| = |H(-k)|$ we find
\begin{equation} \label{1003}
\sum_{\mu X < n \leq X}{\min \left ( 1, \frac{1}{M||n^{c}||}\right )}
		\ll \frac{X\log M}{M} + \frac{\log M}{M}\sum_{1 \leq k \leq M}{|H(k)|} + M\sum_{k > M}{\frac{|H(k)|}{k^{2}}}.
\end{equation}

If $\theta (x) = kx^{c}$, then $\theta '' (x) = c(c - 1)kx^{c - 2} \asymp kX^{c - 2}$ uniformly for
$x \in [\mu X, X]$. Hence, we can apply Van der Corput's theorem (see \cite{17}, chapter 1, Theorem 5) to get
\begin{equation} \label{1004}
H(k) \ll k^{\frac{1}{2}}X^{\frac{c}{2}} + k^{-\frac{1}{2}}X^{1 - \frac{c}{2}}.
\end{equation}
Hence from \eqref{1003} and \eqref{1004} we prove \eqref{5.3}.

When combining Lemma \ref{lemma 10}, \eqref{5.1.1} -- \eqref{5.3} we find that
\begin{align*}
\max _{\Delta \leq \alpha \leq M + 1}{|L(\alpha )|}
\ll x^{\varepsilon }
			\Big (
&
			X^{\frac{1}{3} + \frac{c}{2}}DM^{\frac{1}{2}}
			+ X^{1 - \frac{c}{2}}\Delta ^{-\frac{1}{2}}
			+ X^{\frac{3}{4} + \frac{c}{6}}D^{\frac{2}{3}}M^{\frac{1}{6}} + \\
&
			+ X^{\frac{5}{6}}
			+ X^{1 -\frac{c}{6}}D^{\frac{1}{3}}\Delta ^{-\frac{1}{6}}
			+ X^{1 - \frac{c}{4}}\Delta ^{-\frac{1}{4}}
			+ XM^{-1}
			\Big ).
\end{align*}
Then from last formula, \eqref{3.2} and \eqref{5.1.0} we find
\begin{equation} \label{5.4}
\Gamma _{1}'' \ll
x^{\varepsilon }
			\left (
			X^{\frac{4}{3} + \frac{c}{2} + \delta + \frac{\kappa }{2}}
			+ X^{\frac{7}{4} + \frac{c}{6} + \frac{2\delta }{3} + \frac{\kappa }{6}}
			+ X^{\frac{11}{6}}
			+ X^{2 + \frac{\delta }{3} - \frac{\xi }{6}}
			+ X^{2 - \kappa }
			\right ).
\end{equation}
If we choose $\kappa = \frac{8c - 5}{56}$, then from \eqref{5.1} and \eqref{5.4} we conclude that if
$1 < c < \frac{17}{16}$ then
\begin{equation*}
\Gamma _{1}'' \ll X^{3 - c - \varepsilon }.
\end{equation*}
From \eqref{3.5.1}, \eqref{3.5.2}, \eqref{3.5.3} and \eqref{3.5.7} -- \eqref{3.5.6} we conclude that
\begin{equation} \label{9}
\Gamma \geq |3 \mathcal {N}^{-} - 2\mathcal {N}^{+}|\mathcal ({N}^{+})^{3} B_{1} + O(X^{3 - c}(\log x)^{-4}).
\end{equation}
Now we shall find a lower bound for the difference $3\mathcal {N}^{-} - 2\mathcal {N}^{+}$. It is easy to see that
\begin{equation} \label{10}
\mathcal {B} \asymp (\log X)^{-1}.
\end{equation}
From \eqref{1.02} and \eqref{1.03} we see that
\begin{equation*}
3\mathcal {N}^{-} - 2\mathcal {N}^{+} \geq \mathcal {B}(3f(s_{0}) - F(s_{0})) + O\left ( \log X)^{-\frac{4}{3}} \right ),
\end{equation*}
where $s_{0}$ is defined by \eqref{1.01} and $F(s)$ and $f(s)$ are defined by \eqref{1.04}.
If we choose $s_{0} = 2,95$ then from \eqref{1.01}, \eqref{3.2} and \eqref{5.1} we find
\begin{equation*}
\eta = \frac{\delta }{2,95} = \frac{17 - 16c}{94,4}
\end{equation*}
and also from \eqref{1.04} we find $3f(s_{0}) - F(s_{0}) > 0$.

Now from \eqref{1.02}, \eqref{3.5.7}, \eqref{9} and \eqref{10} we obtain
\begin{equation*}
\Gamma \gg X^{3 - c}(\log X)^{-3}.
\end{equation*}
Therefore $\Gamma > 0$ and this proves Theorem \ref{theorem 1}. \hfill $ \square $ \\

\bibliographystyle{plain}

\vspace{5mm}

\noindent Faculty of Mathematics and Informatics \\
Sofia University "St. Kl. Ohridski" \\
5 J.Bourchier, 1164 Sofia, Bulgaria \\
\vspace{5mm}

\noindent {zhpetrov@fmi.uni-sofia.bg}

\end{document}